\documentclass[a4paper,reqno,12pt]{amsart}
\usepackage{enumerate}
\usepackage{amssymb}
\usepackage{ifthen}
\usepackage{graphicx}
\usepackage{color,xcolor}
\usepackage{mathrsfs}
\usepackage{extarrows}
\newtheorem{thm}{Theorem}
\newtheorem{lem}{Lemma}
\newtheorem{cor}{Corollary}
\newtheorem{defn}{Definition}
\newtheorem{example}{Example}

\newtheorem{claim}{Claim}

\newtheorem{conj}{Conjecture}

\newtheorem{prob}{Problem}

\newenvironment{rem}{%
\bigskip
\noindent \textsl{{\sl Remark. }}}{\bigskip}
\newenvironment{rems}{%
\bigskip
\noindent \textsl{{\sl Remarks. }}}{\bigskip}

\newenvironment{pf}[1][]{%
 \vskip 1mm
 \noindent
 \ifthenelse{\equal{#1}{}}%
  {{\slshape Proof. }}%
  {{\slshape #1.} }%
 }%
{\qed\bigskip}

\newcounter{alphabet}

\newenvironment{Thm}[1][]{\refstepcounter{alphabet}%
\bigskip%
\noindent%
{\bf Theorem \Alph{alphabet}}%
\ifthenelse{\equal{#1}{}}{}{ (#1)}%
{\bf .} \itshape}{\vskip 8pt}

\newenvironment{Lem}[1][]{\refstepcounter{alphabet}%
\bigskip%
\noindent%
{\bf Lemma \Alph{alphabet}}%
{\bf .} \itshape}{\vskip 8pt}

\newcommand{\IC}{{\mathbb C}}
\newcommand{\ID}{{\mathbb D}}

\newcommand{\D}{{\mathbb D}}

\def\be{\begin{equation}}
\def\ee{\end{equation}}
\def\bes{\begin{equation*}}
\def\ees{\end{equation*}}

\newcommand{\bee}{\begin{enumerate}}
\newcommand{\eee}{\end{enumerate}}

\newcommand{\blem}{\begin{lem}}
\newcommand{\elem}{\end{lem}}
\newcommand{\bthm}{\begin{thm}}
\newcommand{\ethm}{\end{thm}}
\newcommand{\bcor}{\begin{cor}}
\newcommand{\ecor}{\end{cor}}
\newcommand{\beg}{\begin{example}}
\newcommand{\eeg}{\end{example}}
\newcommand{\begs}{\begin{examples}}
\newcommand{\eegs}{\end{examples}}
\newcommand{\bdefe}{\begin{defn}}
\newcommand{\edefe}{\end{defn}}
\newcommand{\bprob}{\begin{prob}}
\newcommand{\eprob}{\end{prob}}
\newcommand{\bques}{\begin{ques}}
\newcommand{\eques}{\end{ques}}
\newcommand{\bei}{\begin{itemize}}
\newcommand{\eei}{\end{itemize}}

\newcommand{\bde}{\begin{deter}}
\newcommand{\ede}{\end{deter}}
\newcommand{\bca}{\begin{case}}
\newcommand{\eca}{\end{case}}
\newcommand{\bcl}{\begin{claim}}
\newcommand{\ecl}{\end{claim}}
\newcommand{\bcon}{\begin{conj}}
\newcommand{\econ}{\end{conj}}
\newcommand{\bcons}{\begin{conjs}}
\newcommand{\econs}{\end{conjs}}
\newcommand{\bprop}{\begin{propo}}
\newcommand{\eprop}{\end{propo}}
\newcommand{\br}{\begin{rem}}
\newcommand{\er}{\end{rem}}
\newcommand{\brs}{\begin{rems}}
\newcommand{\ers}{\end{rems}}
\newcommand{\bo}{\begin{obser}}
\newcommand{\eo}{\end{obser}}
\newcommand{\bos}{\begin{obsers}}
\newcommand{\eos}{\end{obsers}}
\newcommand{\bpf}{\begin{pf}}
\newcommand{\epf}{\end{pf}}
\newcommand{\ba}{\begin{array}}
\newcommand{\ea}{\end{array}}
\newcommand{\beq}{\begin{eqnarray}}
\newcommand{\beqq}{\begin{eqnarray*}}
\newcommand{\eeq}{\end{eqnarray}}
\newcommand{\eeqq}{\end{eqnarray*}}

\newcommand{\ds}{\displaystyle}

\newcounter{minutes}
\divide\time by 60
\newcounter{hours}
\multiply\time by 60 \addtocounter{minutes}{-\time}

\begin{document}
\bibliographystyle{amsplain}
\title[Improved Bohr inequality for harmonic mappings]
{Improved Bohr inequality for harmonic mappings}

\def\thefootnote{}
\footnotetext{ \texttt{\tiny File:~\jobname .tex,
          printed: \number\day-\number\month-\number\year,
          \thehours.\ifnum\theminutes<10{0}\fi\theminutes}
} \makeatletter\def\thefootnote{\@arabic\c@footnote}\makeatother

\author[G. Liu]{Gang Liu}
\address{G. Liu, College of Mathematics and Statistics
 (Hunan Provincial Key Laboratory of Intelligent Information Processing and Application),
Hengyang Normal University, Hengyang,  Hunan 421002, China}
\email{liugangmath@sina.cn}

\author[S. Ponnusamy]{Saminathan Ponnusamy
}

\address{S. Ponnusamy, Department of Mathematics,
Indian Institute of Technology Madras, Chennai-600 036, India}
\email{samy@iitm.ac.in}


\keywords{Bohr inequality, Bohr radius, bounded analytic function, harmonic mapping, Schwarz lemma, subordination, quasi-subordination.\\
To appear in Mathematische Nachrichten.
}

\begin{abstract}
Based on improving the classical Bohr inequality, we get in this paper some refined versions for  a quasi-subordination family of functions,
one of which is key to build our results.
By means of these investigations, for a family of harmonic mappings defined in the unit disk $\D$,
we establish an improved Bohr inequality with refined Bohr radius  under particular conditions.
Along the line of extremal problems concerning the refined Bohr radius, we derive a series of results. 
Here the family of harmonic mappings have the form $f=h+\overline{g}$,
where $g(0)=0$, the analytic part $h$ is bounded by 1
and that $|g'(z)|\leq k|h'(z)|$ in $\D$ and for some $k\in[0,1]$.
\end{abstract}

\maketitle \pagestyle{myheadings}
\markboth{G. Liu and S. Ponnusamy}{Improved Bohr inequality for harmonic mappings}

\section{Introduction}  \label{sec1}
Throughout the paper, ${\mathcal B}$ denotes the set of all analytic functions $f$ in the unit disk  $\mathbb{D}=\{z\in \IC:\, |z| < 1 \}$
such that $|f(z)|\leq1$ for all $z\in\D$. As with the standard decomposition of complex-valued harmonic functions (cf. \cite{CS,dur2004,PR}),
let $\mathcal{H}$ and $\mathcal{H}_k$ denote the set of harmonic mappings defined by
$$\mathcal{H}=\left \{f=h+\overline{g}:\, h~\text{and}~g ~\text{are ~analytic ~in} ~\D ~\text{with}~ g(0)=0\right\}
$$
and
$$\mathcal{H}_k=\left \{f=h+\overline{g} \in\mathcal{H}:\, h\in{\mathcal B}~\text{and}~|g'|\leq k|h'| ~\text{in}~ \D ~\text{for some}~ k\in[0,1]\right\},
$$
respectively. Clearly, $\mathcal{H}_0\equiv {\mathcal B}$.
Let us recall few basic notions about harmonic mappings.
A function $f=h+\overline{g} \in \mathcal{H}$ is sense-preserving whenever $J_f=|h'|^2-|g'|^2>0$ in $\ID$, or
equivalently $h'(z) \neq 0$ and $|g'(z)|<|h'(z)|$ for all $z\in \mathbb{D}$.
Further, if its dilatation $\omega_f=g'/h'$ satisfies $|\omega_f|\leq k<1$ in $\mathbb{D}$,
then $f$ is called a $K-$quasiregular mapping, where $K=(1+k)/(1-k)$.
For more details of the importance, background, development and results, we refer to the monograph of
Duren \cite{dur2004} and the survey article of Ponnusamy and Rasila \cite{PR}.

Let us recall the classical theorem of Herold Bohr.

\begin{Thm}\label{thmA}
{\rm (\cite[Bohr (1914)]{boh})}
Suppose that $f\in {\mathcal B}$ and $f(z)=\sum_{n=0}^{\infty} a_n z^n$.
Then the following sharp inequality holds:
$$\sum_{n=0}^{\infty} |a_n|\, r^n \leq 1
~\mbox{ for  $r\leq 1/3$.}
$$
\end{Thm}

In recent years, a number of researchers revisited the work of Bohr--improving and extending this work to more general setting.
Bohr's original proof had the above mentioned inequality only for $r\leq1/6$, which was later improved
independently by  M. Riesz, I. Schur, F. Wiener and some others. 
We call the sharp constant $1/3$ in Theorem~A 
the Bohr radius for the family ${\mathcal B}$. Later
proofs were given by Sidon \cite{sid} and Tomic \cite{tom}. See also \cite{PPS,PS2004,PS2006} and
and the recent survey chapters \cite{AAP} and \cite[Chapter 8]{GarMasRoss-2018}.
In addition, if $|a_0|$ in Bohr inequality is replaced by $|a_0|^p$, where $1\leq p\leq 2$
then the constant $1/3$ could be replaced by $p/(2+p)$, see \cite[Proposition~1.4]{bla}.
In \cite[Remark 1]{PVW2}, this result was shown to be true in refined form even for the extended range $0<p\leq 2$.
Moreover, if $a_0=0$ in Theorem~A 
then the sharp Bohr radius is improved to be $1/\sqrt{2}$ which was shown by
Bombieri \cite{Bomb-1962} in 1962. See also \cite{KP2017,KayPon_AAA18}, \cite[Corollary~2.9]{PPS} and the recent paper of Ponnusamy and Wirths \cite{PW}
where one can find this result as a special case of each of theirs.

It is astonishing that various generalizations of the classical Bohr inequality
have been investigated in different branches of mathematics.
For instance, Hardy spaces \cite{BDK5,DjaRaman-2000}, Bloch spaces \cite{KPS,LP}, harmonic mappings \cite{abu,AANgH,AKP2019,EPR,KP2018-1,KPS,LPW,LiuP},
Dirichlet series \cite{BaluCQ-2006},  logarithmic power series \cite{BhowDas-19}, functions in Banach space \cite{bla}, and holomorphic
functions of several variables  \cite{Aizen-00-1,aiz,Aiz07,BDK5,BoasKhavin-97-4,DF,DFOOS,DjaRaman-2000}.

To prove or improve the classical Bohr inequality, one mainly relies on the sharp coefficient inequalities. In fact, Theorem~A 
can be easily deduced from the classical result $|a_n|\leq1-|a_0|^2$ $(n\geq1,~f\in{\mathcal B})$.
As mentioned for example in \cite{PW}, this inequality follows quickly from a result on subordination due to Rogosinski.
However,  its sharpness cannot be obtained in the extremal case $|a_0|<1$, which was pointed out in \cite{KP2017}.
Therefore, on one hand, the sharp version of Theorem~A 
has been achieved for any individual function from
${\mathcal B}$ (see \cite{AKP} and some subclass of univalent functions (see \cite{abu,AANgH})).
On the other hand, through a refined version of the coefficient inequalities found by Carlson (see \cite{Car}),
 Bohr's inequality was refined and improved  in the following way (see also \cite{PVW2}).

In what follows we let $\|f\|_r^2=\sum_{n=1}^\infty|a_n|^2r^{2n}$ whenever
$f(z)=\sum_{n=1}^{\infty} a_n z^n$ converges for $|z|<1$ and $r<1$.

\begin{Thm}{ \rm(\cite[Theorem~2]{PVW})} \label{thmB}
 Suppose that $f\in {\mathcal B}$,  $f(z)=\sum_{n=0}^{\infty} a_n z^n$  and $f_0(z)=f(z)-f(0)$.
Then for $p=1,2$, we have the following sharp inequality
\begin{equation*}
|a_0|^p+ \sum_{n=1}^\infty |a_n|r^n+  \frac{1}{1+|a_0|}\left(\frac{1+|a_0|r}{1-r}\right)\|f_0\|_r^2\leq 1 ~\mbox{for}~ r \leq \frac{1}{1+(1+|a_0|)^{2-p}}.
\end{equation*}
\end{Thm}

Besides these results, there are a number of works about Bohr inequality for the family ${\mathcal B}$.
One is to consider Bohr inequality for functions of the form $f_{p,m}(z)=\sum_{n=0}^{\infty} a_{pn+m} z^{pn+m}$ in ${\mathcal B}$
(see \cite{ABS,KP2017,KP2018-1,LLS}).
In particular, $f_{p,1}$ is called $p-$symmetric function and $f_{2,1}$ is called odd function.
The other is to study the Bohr-Rogosinski inequality (see \cite{AKP2019,KP-pre,LLS}), which was introduced by Kayumov and Ponnusamy
in \cite{KP-pre} based on the notion of Rogosinski inequality investigated in \cite{LG,rog,SS}.
Another aspect of it is to build  different Bohr type inequalities associated with alternating series,
area, modulus of $f$ or $f-a_0(f)$,  and higher order derivatives of $f$ in part or in whole etc.. These include the
works of  \cite{ABS,AKP2019,IKP-20,KP2018-2,LLS,LSX}. For some other related topics may be found in  \cite{LPW,LiuP}.
As mentioned above, there exist Bohr's theorems to more general domains or higher dimensional spaces, holomorphic functions defined
on bounded complete Reinhardt domain in $\IC^n$, and operator-theoretic Bohr radius. See for example, \cite{Aizen-00-1,aiz,Aiz07,BoasKhavin-97-4,HH2012}.

There are few harmonic extensions concerning Bohr inequality for the family ${\mathcal H}_k$.
It was first considered in the paper \cite{KPS} for ${\mathcal H}_k$ $(k\neq1)$
and a couple of problems on Bohr's inequality for its subclass were also posed. Here, it should be mentioned that this work was
motivated by the work from \cite{abu, AANgH}.
The problems proposed in \cite{KPS} were solved in \cite{BD,LPW} (see also \cite{LiuP}) by means of quasi-subordination with special forms,
which was generalized  in order to get more results of Bohr inequality for ${\mathcal H}_k$ in \cite{AKP}.
It is emphasized that their proofs are heavily depending on Theorem~A. 
Let us now recall the following.

\bdefe
For any two analytic functions $f$ and $g$ in $\mathbb{D}$, we say that the function
$f$ is quasi-subordinate to $g$ (relative to $\Phi$),
denoted by $f(z) \prec_{q} g(z)$ in $\mathbb{D}$ if there exist two functions
$\Phi\in\mathcal{B},~\omega\in\mathcal{B}$ with $\omega(0)=0$ such that $f(z)=\Phi(z) g(\omega(z))$.
\edefe

There are two special cases which are of particular interest.
The choice $\Phi(z)=1$  corresponds to subordination which is denoted by $f \prec g$,
whereas $\omega(z)=z$ gives majorization, i.e. reduces to the form $f(z)=\Phi(z) g(z)$, which is equivalent to $|f(z)|\leq|g(z)|$ in $\D$.
Note that $g'$ is majorized by $kh'$ in the definition of ${\mathcal H}_k$.
Along the lines of works on Bohr inequality for the family ${\mathcal B}$ in \cite{KP2018-2},
few different formulations of improved Bohr inequalities for ${\mathcal H}_1$ and ${\mathcal H}_k$
were obtained in \cite{EPR} and \cite{AKP2019}, respectively.
For more recent advances on Bohr's inequality for the family ${\mathcal H}$, the reader may refer for example,
\cite{KKP-2020,KP2018-1,K-Sahroo-pre,LP,LmsP,LPW,LiuP}.

Now, a variety of Bohr radii exist because of different formulations and refinements (cf. \cite{KP2018-1,KP2018-2},
and \cite[Theorem~2]{PVW}, i.e., Theorem~B) 
of the classical Bohr
inequality, and thus it becomes more and more complex in some situation as you see in our investigation in this paper, especially when we try to study
the extension of Bohr inequality from ${\mathcal B}$ to ${\mathcal H}_k$.

For the family ${\mathcal B}$, we know that the classical Bohr radius is a constant,
which is improved to be a function of the modulus of the constant term (see Theorem~B). 
For the family ${\mathcal H}_k$, the expression of sharp Bohr radius either is a constant or a function of the variable $k$.
Furthermore, it is worth pointing out that there is only one result  related to both the constant $k$
and the modulus of the constant term of its analytic part, but such result holds with additional assumptions
(cf. \cite[Theorem~2.9]{AKP}).

In view of these reasonings, some interesting questions emerge.
In the process of harmonic extension, it is natural to ask whether the formulation of Bohr inequality is complex
so that it can cover or improve some known results or not? Equivalently, we ask under what conditions, the Bohr radius will be depending
on $k$ or the modulus of the constant term of its analytic part, or both?
Another natural question is to improve Bohr inequality or Bohr radius, and to integrate some of the known results into simplified forms.
In this paper, we try to answer these questions partly.

The paper is organized as follows. In Section \ref{sec2}, we improve the classical Bohr inequality and obtain some refined versions for a
quasi-subordination family of functions in Section \ref{sec3}. In view of these investigations, improved Bohr type inequalities for ${\mathcal H}_k$ are
established in Section \ref{sec4}.

The proofs of our results rely on a couple of lemmas which we recall now.

\begin{Lem} { \rm(\cite[Proof~ of~ Theorem ~1]{KP2017} and \cite{KP2018-1})} \label{lem1}
 Suppose that $f\in {\mathcal B}$ and $f(z)=\sum_{n=0}^{\infty} a_n z^n$.
 Then we have
$$ \sum_{n=1}^{\infty}\left|a_{n}\right| r^{n} \leq
\begin{cases}
\displaystyle A(r):=r\frac{1-|a_0|^2}{1-r|a_0|} & \text{for}\quad |a_0|\geq r,\\[4mm]
 \displaystyle B(r):=r\frac{\sqrt{1-|a_0|^2}}{\sqrt{1-r^2}} &\ \text{for}\quad |a_0|<r.
\end{cases}
$$
\end{Lem}

\begin{Lem}\label{lem2} { \rm (\cite[p.107, Proof of Theorem~1]{PVW})}
Suppose that $f\in {\mathcal B}$, $f(z)=\sum_{n=0}^{\infty} a_n z^n$ and $f_0(z)=f(z)-f(0)$.
Then we have
$$	\sum_{n=1}^{\infty}|a_n|r^n+ \frac{1}{1+|a_0|}\left(\frac{1+|a_0|r}{1-r}\right)\|f_0\|_r^2\leq(1-|a_0|^2)\frac{r}{1-r}\quad \text{for}~ r\in[0,1).
$$
\end{Lem}

\begin{Lem}\label{lem3}
{\rm (Schwarz-Pick lemma)}
Suppose that $f\in {\mathcal B}$. Then we have
$$|f(z)|\leq\frac{|z|+|f(0)|}{1+|f(0)|\, |z|}\quad \text {and}\quad
 |f'(z)|\leq\frac{1-|f(z)|^2}{1-|z|^2}\quad \text{for } z\in\D.
$$
Equality holds at some point $z_0\in \ID$ either in the first inequality or in the second inequality, if and only if
 $f(z)=c\frac{z+a}{1+\overline{a}\,z}$, $z\in \ID$, for some $c$ with $|c|=1$ and $a\in \ID$.
\end{Lem}

\section{Improved versions of  the classical Bohr inequality }  \label{sec2}

In what follows, for the sake of simplicity, we denote three functions as following:
$$\omega_a(z)=\frac{z+a}{1+az}= a +(1-a^2)\sum_{k=1}^\infty (-a)^{k-1}z^k,
\quad z\in\D,~a\in[0,1),
$$
$$t_p(x)=\frac{1-x^2}{1-x^p},\quad x\in[0,1),
$$
and
$$ r_p(x)=
\begin{cases}
\displaystyle \frac{1-x^p}{\sqrt{1-x^2+(1-x^p)^2}}=\frac{1}{\sqrt{1+\frac{1}{1-x^2}t_p^2(x)}} & \text{for}\quad x\in[0,C(p)),\\[4mm]
 \displaystyle \frac{1-x^p}{1-x^2+x(1-x^p)}=\frac{1}{x+t_p(x)} &  \text{for}\quad x\in[C(p),1),\\[4mm]
  \displaystyle \frac{p}{2+p} & \text{for}\quad x=1,
\end{cases}
$$
where $p>0$ and $C(p)$ is the unique solution of the equation $1-x-x^p=0$ in the interval $(0,1)$.
Clearly, $r_p(0)=1/\sqrt{2}$ for all $p>0$. We observe that
$$r_p(x)\leq1/\left(x+\frac{1-x^2}{x}\right)=x<1,\quad \text{for~}x\in[C(p),1),
$$
which implies $r_p(x)<1$ for $x\in[0,1]$. Clearly, $C(1)=1/2$ and
\be \label{LP8-eq1}
r_1(x)=
\begin{cases}
\displaystyle \sqrt{\frac{1-x}{2}} & \text{for}\quad x\in[0,\frac{1}{2}),\\[4mm]
 \displaystyle \frac{1}{1+2x} & \text{for}\quad x\in[\frac{1}{2},1].
\end{cases}
\ee
The following properties of $t_p$ and $r_p$ will be always used later and we leave it as an exercise.

\blem\label{lem4} For the functions $t_p$ and $r_p$ defined as above, we have the following:
\begin{enumerate}
\item[{\rm (a)}]
The function $t_p$ {\rm(}resp. ~$r_p${\rm)} is continuous in the interval $[0,1)$ {\rm(}resp. $[0,1]${\rm)}.

\item[{\rm (b)}] For each $p\in(0,2)$ \rm{(}resp. $p>2${\rm)}, the function $t_p$
is strictly increasing \rm{(}resp. decreasing{\rm)} in $[0,1)$ and $t_p\in[1,2/p)$ \rm{(}resp. $t_p\in(2/p,1]${\rm)}.
\end{enumerate}
In particular, the function $r_p$ is strictly decreasing from $1/\sqrt{2}$ to $p/(2+p)$ in $[0,1]$ when $p\in(0,2]$.
\elem

\bthm \label{thm1}
Suppose that $p>0$ and $f\in {\mathcal B}$ with $f(z)=\sum_{k=0}^{\infty} a_k z^k$.
Then $$D_f^p(z):=|a_0|^p+\sum_{k=1}^{\infty}|a_k|r^k \leq 1\quad \text{for} \quad r=|z|\leq r_p(|a_0|),
$$
and $r_p(|a_0|)$ cannot be improved for each $p>0$ if $|a_0|\in[C(p),1)\cup\{0\}$.
\ethm
\bpf
Fix $p>0$ and set $a=|a_0|$. Clearly, $a\leq1$.
The proof is trivial if $a=1$, since $f(z)=ae^{i\theta}$ for some $\theta\in\mathbb{R}$.
We only consider the case of $a\in[0,1)$.
Note that $r_p(a)\leq a$ when $a\in[C(p),1)$.
It follows from Lemma~C 
that
$$D_f^p(z) \leq a^p+A(r)\leq a^p+A\left(r_p(a)\right)=1~ \mbox{ for $r\leq r_p(a)$ and $a\in[C(p),1)$.}
$$
For $a\in[0,C(p))$, we observe that $1-a^p>a$ so that $(1-a^p)^2>a^2$, which means that $a<r_p(a)$.
It follows from Lemma~C 
again that
$$D_f^p(z) \leq a^p+A(r)\leq a^p+A(a)=a^p+a<1 ~ \mbox{ for $r\leq a$ and $a\in[0,C(p))$,}
$$
and
$$D_f^p(z) \leq a^p+B(r)\leq a^p+B(r_p(a)) = 1~ \mbox{ for $a< r\leq r_p(a)$ and $a\in[0,C(p))$.}
$$

It remains to show the sharpness part.
If $|a_0|\in[C(p),1)$,  then the extremal function can be chosen as $\omega_a$ with $a\in[C(p),1)$.
For this function,  simple computations show that
$$D_{\omega_a}^p(z)=a^p+(1-a^2)\sum_{k=1}^\infty a^{k-1}r^k=a^p+\frac{(1-a^2)r}{1-ar},
$$
which is bigger than $1$ is equivalent to the condition $r>r_p(a)$.

If $a_0=0$, then we consider the function  $f(z)=z\omega_b(z)$ with $b=1/\sqrt{2}$
and obtain by elementary calculations that
$$D_f^p(z)=br+(1-b^2)\sum_{k=1}^\infty b^{k-1}r^{k+1}=br+\frac{(1-b^2)r^2}{1-br}=\frac{br}{1-br},
$$
which is bigger than $1$ is equivalent to the condition $r>1/\sqrt{2}$.
This completes the proof of the theorem.
\epf

\bcor \label{cor1}
{\rm (See \cite[Remark 1]{PVW2} in refined form)}
Suppose that $p\in(0,2]$ and $f\in {\mathcal B}$ with $f(z)=\sum_{k=0}^{\infty} a_k z^k$.
Then the following sharp inequality holds:
$$|a_0|^p+\sum_{k=1}^{\infty}|a_k|r^k \leq 1\quad \text{for} \quad r\leq r_p(1)=\frac{p}{2+p}.
$$
\ecor

 We would like to point out that Corollary \ref{cor1} was obtained in \cite[Proposition~1.4]{bla} for $p\in[1,2]$,
which was generalized to the case $0<p\leq2$ in a refined form in \cite{PVW2}.
Moreover, the constant $p/(2+p)$ in Corollary \ref{cor1} is the minimum of the function $r_p(x)$ in the interval $[0,1]$,
which is difficult to compute in the case $p>2$.
In fact, the monotonicity of $r_p$ is very complex when $p>2$. For instance, simple computations show that
$$r_4(1/2)>r_4(1/3)>r_4(0)>r_4(1)>1/2.
$$

\section{Refined versions for a quasi-subordinating family of functions}  \label{sec3}

In this section, on the basis of Theorem \ref{thm1}, we obtain a refined version of \cite[Theorem~2.1]{AKP} for a quasi-subordinating
family of functions. In order to present its proof, we need precise relationships concerning quasi-subordination and this is done using
the approach of  \cite[Proof of Theorem~2.1]{AKP}. Moreover, for the proof of Theorem  \ref{thm2} and its corollaries that follow, it has become
necessary to indicate the major steps in brief.

\bthm \label{thm2}
 Let $f(z)$ and $g(z)$ be two analytic functions in $\mathbb{D}$ with the Taylor series expansions $f(z)=\sum_{k=0}^{\infty} a_{k} z^{k}$ and $g(z)=\sum_{k=0}^{\infty} b_{k} z^{k}$.
If there exist two analytic functions $\Phi\in\mathcal{B}$ and $\omega\in\mathcal{B}$ with $\omega(0)=0$ 
such that $f(z)=\Phi(z) g(\omega(z))$ in $\D$.
Then
\[
\sum_{k=0}^{\infty}\left|a_{k}\right| r^{k} \leq \sum_{k=0}^{\infty}\left|b_{k}\right| r^{k}
\quad \text {for} \quad r \leq \min\{r_1(|\Phi(0)|),~r_1(|\omega'(0)|)\},
\]
where $r_1(x)$ is defined by \eqref{LP8-eq1}. 
\ethm
\bpf
Let $\omega(z)=\sum_{n=1}^{\infty} \alpha_{n} z^{n}$. Then, for $k \in \mathbb{N},$ we can write
\[
\omega^{k}(z)=\sum_{n=k}^{\infty} \alpha_{n}^{(k)} z^{n}=z^k(\alpha_1^k+\cdots).
\]
Since $\omega\in\mathcal{B}$ with $\omega(0)=0$,  we have $|\alpha_1|^k\leq|\alpha_1|=|\omega'(0)|\leq1$ for all $k\in\mathbb{N}$.
It follows from Theorem \ref{thm1} that
\[
\sum_{n=k}^{\infty}\left|\alpha_{n}^{(k)}\right| r^{n-k} \leq 1 \text { for } r \leq r_1(|\alpha_1|^k)
\text{ and } k\in\mathbb{N},
\]
and, because $r_1(x)$ is decreasing monotonically in $[0,1)$ by Lemma \ref{lem4}, this implies
\be\label{equ1}
\sum_{n=k}^{\infty}\left|\alpha_{n}^{(k)}\right| r^{n-k} \leq 1 \text { for } r \leq r_1(|\alpha_1|) \text{ and } k\in\mathbb{N}.
\ee
Writing $\Phi(z)=\sum_{m=0}^{\infty} \phi_{m} z^{m}$, by Theorem \ref{thm1}, we have
\be\label{equ2}
\sum_{m=0}^{\infty}\left|\phi_{m}\right| r^{m} \leq 1 \text { for } r \leq r_1(|\Phi(0)|).
\ee
For simplicity, we introduce  $\omega^{0}(z)=1=\sum_{n=0}^{\infty} \alpha_{n}^{(0)} z^{n}$, where $\alpha_{0}^{(0)}=1$, $\alpha_{n}^{(0)}=0$ for $n \geq 1$.
Then, as in \cite[Theorem~2.1]{AKP}, we can rewrite the relation $f(z)=\Phi(z)g(\omega(z))$ equivalently in terms of power series as
$$\sum_{k=0}^{\infty} a_{k} z^{k} = \sum_{k=0}^{\infty}\left(\sum_{m+j=k} \phi_{m} B_{j}\right) z^{k}, \quad a_{k}=\sum_{m+j=k} \phi_{m} B_{j} \text { for each } k \geq 0,
$$
where $B_{k}=\sum_{n=0}^{k} b_{n} \alpha_{k}^{(n)}.$ 
Applying the triangle inequality, we easily have
\[
\begin{aligned}
\sum_{k=0}^{\infty}\left|a_{k}\right| r^{k}  & \leq
\left(\sum_{m=0}^{\infty}\left|\phi_{m}\right| r^{m}\right) \sum_{k=0}^{\infty}\left|B_{k}\right| r^{k} \\
& \leq \sum_{k=0}^{\infty}\left|B_{k}\right| r^{k} \text { for } r \leq r_1(|\Phi(0)|), \quad \mbox{(by \eqref{equ2})}.
\end{aligned}
\]
Also, because $\left|B_{k}\right| \leq \sum_{n=0}^{k}\left|b_{n}\right|\big |\alpha_{k}^{(n)}\big|,$ we obtain that
\[
\begin{aligned}
\sum_{k=0}^{\infty}\left|B_{k}\right| r^{k} & \leq \sum_{k=0}^{\infty} \sum_{n=0}^{k}\left|b_{n}\right|\big |\alpha_{k}^{(n)}\big| r^{k}=\sum_{k=0}^{\infty}\left|b_{k}\right| \sum_{n=k}^{\infty}\left|\alpha_{n}^{(k)}\right| r^{n} \\
&=\sum_{k=0}^{\infty}\left|b_{k}\right|\left(\sum_{n=k}^{\infty}\left|\alpha_{n}^{(k)}\right| r^{n-k}\right) r^{k} \\
& \leq \sum_{k=0}^{\infty}\left|b_{k}\right| r^{k} \text { for } r \leq r_1(|\omega'(0)|),\quad \mbox{(by \eqref{equ1})},
\end{aligned}
\]
and hence,
\[
\sum_{k=0}^{\infty}\left|a_{k}\right| r^{k} \leq \sum_{k=0}^{\infty}\left|B_{k}\right| r^{k} \leq \sum_{k=0}^{\infty}\left|b_{k}\right| r^{k} \text { for } r \leq \min\{r_1(|\Phi(0)|),r_1(|\omega'(0)|)\}.
\]
The proof of Theorem \ref{thm2} is complete.
\epf

\begin{cor} \label{cor2}
Suppose that $f\prec g$, where $f$ and $g$  
are defined as in Theorem {\rm \ref{thm2}}. Then we have
\begin{enumerate}
	\item[{\rm (a)}] $\ds \sum_{k=0}^{\infty}\left|a_{k}\right| r^{k} \leq \sum_{k=0}^{\infty}\left|b_{k}\right| r^{k}$ $\ds
\text{ for } r \leq r_1(|a_1/b_1|),$
when $b_1\neq0$.
	\item[{\rm (b)}] $\ds \sum_{k=0}^{\infty}\left|a_{k}\right| r^{k} \leq \sum_{k=0}^{\infty}\left|b_{k}\right| r^{k} $ $\ds
\text{ for } r \leq 1/3,$
when $b_1=0$.
\end{enumerate}
Moreover, $r_1(|a_1/b_1|)$ cannot be improved if  $|a_1/b_1|\in[1/2,1)\cup\{0\}$, and the constant $1/3$
 in {\rm (b)} cannot be improved.
\end{cor}

\bpf
\rm{(a)} Let $b_1\neq0$ and  $f\prec g$. 
Then  $f(z)=\Phi(z)g(\omega(z))$ in $\D$,
where $\Phi(z)=1$ and $\omega\in\mathcal{B}$ with $\omega(0)=0$.
Now $f'(z)=g'(\omega(z))\omega'(z)$, which implies that $\omega'(0)=f'(0)/g'(0)=a_1/b_1$.
As $\Phi(z)=1$, from the proof  of Theorem \ref{thm2},
the desired result follows with the replacement of  $r\leq  r_1(|\Phi(0)|)$ in \eqref{equ2} by $r<1$.

For the sharpness part, we first let $|a_1/b_1|\in[1/2,1)$ and consider $g(z)=b_1z$ and
$f(z)=b_1z\omega_a(z)$ with $a=|a_1|/|b_1|$.
Then $f \prec g$  and it is easy to see that
$$\sum_{k=0}^{\infty}\left|a_{k}\right| r^{k} \leq \sum_{k=0}^{\infty}\left|b_{k}\right| r^{k}
\quad \text{if and only if}\quad r\leq r_1(a)=\frac{1}{1+2a}.
$$

Next we let $a_1=0$. In this case choose $g(z)=b_1z$ with $b_1\neq0$ and $f(z)=b_1z^{2}\omega_a(z)$ with $a=1/\sqrt{2}$.
Again $f \prec g$ and it is easy to see that
$$\sum_{k=0}^{\infty}\left|a_{k}\right| r^{k} \leq \sum_{k=0}^{\infty}\left|b_{k}\right| r^{k}
\quad \text{if and only if} \quad r\leq r_1(0)=\frac{1}{\sqrt{2}}.
$$

\rm{(b)}
Note that $b_1=0$ and $g(z)\not\equiv0$. Again $f(z)= g(\omega(z))$ in $\D$ for some $\omega\in\mathcal{B}$ with $\omega(0)=0$.
The result follows from Theorem \ref{thm2} with $\Phi(z)=1$, since $r_1(x)\geq1/3$ for $x\in[0,1]$.

Next we will show the part of sharpness. Let us consider function
$$f(z)=z^2\left( \frac{z-a}{1-az}\right )^2 =z^2\sum_{k=0}^\infty A_kz^k,
$$
where $a\in(0,1)$, $A_0=a^2$ and $A_k= (1-a^2) a^{k-2}(k-1-(k+1)a^2)  $. Then $f(z)\prec z^2$ in $\D$.

To search for the upper bound of $r$ in the inequality $r^2\sum_{k=0}^{\infty}\left|A_{k}\right| r^{k} \leq r^2$,
it suffices to consider that of $r$ in the inequality $\sum_{k=0}^{\infty}\left|A_{k}\right| r^{k} \leq 1$.
We observe that if $\frac{N-1}{N+1}\leq a^2<\frac{N}{N+2}$ 
for some $N\in\mathbb{N}$, then $A_k\leq0$ for $k\leq N$, and
$A_k>0$ for $k> N$, and hence,  we can write
$$ S_{a,N}(r)  =\sum_{k=0}^{\infty}\left|A_{k}\right| r^{k} =a^2-\sum_{k=1}^N A_kr^k+\sum_{k=N+1}^\infty A_kr^k=\left(\frac{r-a}{1-ar}\right)^2-2\sum_{k=1}^N A_kr^k.
$$
We denote the upper bound of $r$ in the inequality $S_{a,N}(r)\leq1$ by $r(a,N)$.
Next we will show
$$\inf_{a\in(1/\sqrt{2},1)}r(a,N)=1/3.
$$
It follows from Theorem \ref{thm2} that $r(a,N)\geq r_1(a^2)\geq1/3$.
Note that $N$ increases to $+\infty$ when $a$ approaches  $1$.
To certify our assertion, we introduce
\[
\begin{aligned}
S_a(r)&=\left(\frac{r-a}{1-ar}\right)^2-2\sum_{k=1}^\infty A_kr^k\\
&=\left(\frac{r-a}{1-ar}\right)^2+2(1-a^2)\sum_{k=1}^\infty a^{k-2}((k+1)a^2-(k-1))r^k.
\end{aligned}
\]
By computation,  we get that for $a>1/\sqrt{2}$,
\[
\begin{aligned}
S_a(r) &=\left(\frac{r-a}{1-ar}\right)^2+2(1-a^2)\left(\frac{ar(2-ar)}{(1-ar)^2}-\frac{r^2}{(1-ar)^2}\right)\\
&=1-\frac{1-a^2}{(1-ar)^2}((1+2a^2)r^2-4ar+1)\\
&=1-\frac{(1-a^2)(1+2a^2)}{(1-ar)^2}(r-\alpha_{+})(r-\alpha_{-}),\quad \alpha_{\pm} =\frac{1}{2a\pm \sqrt{2a^2-1}}.
\end{aligned}
\]
In the above sum, we have used the formula
$$\sum_{k=N}^\infty kz^{k-1}=\frac{z^{N-1}}{(1-z)^2}(N+(1-N)z) \quad \text {for~} N\in\mathbb{N}.
$$
It is easy to see that $S_a(r)\leq1$ if and only if $r\leq\alpha_+$ or $\alpha_{-}\leq r<1$.
We observe that $S_{a,N}(r)\geq S_a(r)$ for all $a\in(1/\sqrt{2},1)$, and thus, $r(a,N)\leq\alpha_+$ for $a\in(1/\sqrt{2},1)$.
We find that
$$\inf_{a\in(1/\sqrt{2},1)}\alpha_+=\frac{1}{3}
$$
which yields
$$ \inf_{a\in(1/\sqrt{2},1)}r(a,N)=\frac{1}{3}.
$$
The proof of Corollary \ref{cor2} is finished.
\epf

\begin{cor} \label{cor3}
Suppose that $|f(z)|\leq|g(z)|$ for all $z \in \mathbb{D}$, where $f$ and $g$
are defined as in Theorem {\rm \ref{thm2}} with $b_k\neq0$ for some non-negative integer $k$. Then we have
\[
\sum_{k=0}^{\infty}\left|a_{k}\right| r^{k} \leq \sum_{k=0}^{\infty}\left|b_{k}\right| r^{k}
\text { for } \quad r \leq r_1(|a_q/b_q|),
\]
where $q$ is the order of the zero of $g$ at $0$. Moreover, $r_1(|a_q/b_q|)$ cannot be improved if $|a_q/b_q|\in[1/2,1)\cup\{0\}$.
\end{cor}
\bpf
Suppose that $|f(z)|\leq|g(z)|$ for all $z \in \mathbb{D}$. Then $f$ can be written as
$f(z)=\Phi(z)g(\omega(z))$ in $\D$,
where  $\omega (z)=z$ and $\Phi =f/g$ is an analytic function with $|\Phi(z)|\leq1$ in $\D$, and
$\Phi(0)=a_q/b_q$. By the method of the proof of Theorem \ref{thm2}, the desired result follows with the
replacement of  $r\leq  r_1(|\omega'(0)|)$ in \eqref{equ1} by $r<1$.

For the sharpness part, if $|a_q/b_q|\in[1/2,1)$, then we consider $g(z)=b_qz^q$ and
$f(z)=b_qz^q\omega_a(z)$ with $a=|a_q|/|b_q|$,
where $q$ is a non-negative integer.
Clearly $|f(z)|\leq|g(z)|$ for all $z\in\D$, and it is easy to see that
$$\sum_{k=0}^{\infty}\left|a_{k}\right| r^{k} \leq \sum_{k=0}^{\infty}\left|b_{k}\right| r^{k}
\quad \text{if and only if}\quad r\leq r_1(|\Phi(0)|)= r_1(a)=\frac{1}{1+2a}.
$$
If $\Phi(0)=0$, i.e., if $a_q=0$, then we choose $g(z)=b_qz^q$ with $b_q\neq0$ and
$f(z)=b_qz^{q+1}\omega_a(z)$ with $a=1/\sqrt{2}$.
It is easy to see that $|f(z)|\leq|g(z)|$ for all $z\in\D$ and
$$\sum_{k=0}^{\infty}\left|a_{k}\right| r^{k} \leq \sum_{k=0}^{\infty}\left|b_{k}\right| r^{k}
\quad \text{if and only if} \quad r\leq r_1(0)=\frac{1}{\sqrt{2}}.
$$
The proof of Corollary \ref{cor3} is finished.
\epf

\brs
Recall that $r_1(x)\geq1/3$ for all $x\in[0,1]$.
Thus, Theorem \ref{thm2}, Corollary \ref{cor2} and Corollary \ref{cor3} are refined versions of Theorem 2.1,
Corollary 2.2 (\cite[Lemma~1]{BD}) and Corollary 2.3 in \cite{AKP}, respectively.
\ers

\section{Improved  Bohr inequality for harmonic mappings}  \label{sec4}

In this section, in order to improve some Bohr inequalities for the family $\mathcal{H}_k$,
we need a key lemma. For simplicity, we introduce some notations. Suppose that
\beq \label{equ3}
f(z)=h(z)+\overline{g(z)}=\sum_{n=0}^{\infty} a_{n} z^{n}+\overline{\sum_{n=1}^{\infty} b_{n} z^{n}}
\eeq
is a harmonic mapping in $\mathbb{D}$. Without special statement,  let $h$ be not identically a constant and
$h_0(z)=h(z)-h(0)$. We define the quantity $E_f(k,r)$  by
$$E_f(k,r)=\sum_{n=1}^{\infty}\left|a_{n}\right| r^{n}+\sum_{n=1}^{\infty}\left|b_{n}\right| r^{n}
+ \frac{1}{1+|a_0|}\left(\frac{1+|a_0|r}{1-r}\right)\left(||h_0||_r^2+c(k)||g||_r^2\right),
$$
where $k\in[0,1]$, $r=|z|\in[0,1)$
and
$$ c(k)=
\begin{cases}
\displaystyle 0 & \text{for }~ k=0,\\
 \displaystyle 1/k & \text{for }~ k\in(0,1].
\end{cases}
$$

\blem\label{lem5}
Suppose that $f=h+\overline{g}$ is a harmonic mapping in $\mathbb{D}$ with the form of \eqref{equ3} such that $f\in\mathcal{H}_k$.
Then for   $k\in (0,1]$, we have
$$E_f(k,r) \leq (1-|a_0|^2)(1+k)\frac{r}{1-r},\quad
\text{ for } r=|z|\leq r_1(|b_q|/(k|a_q|)),
$$
where $q$ is the order of the zero of $h_0$ at $0$. For $k=0$, the above inequality holds for $r<1$.
Furthermore, equality holds  in the above inequality if $f=h+\overline{g}$ with $h=\omega_a$ and $g=k(\omega_a-a)$.
\elem
\bpf
For $k=0$,  the conclusion is a direct consequence of Lemma ~D. 
Now we only consider the case of $k\in(0,1]$.
If $q$ is the order of zero of $h_0$ at $0$, then  $h$ can be rewritten as
$$h(z)=h(0)+\sum_{n=q}^\infty a_nz^n \quad (a_q\neq0).
$$
Clearly, $kh'(z)=k\sum_{n=q}^\infty na_nz^{n-1}$.
Since $|g^{\prime}(z)| \leq k|h^{\prime}(z)|$ in $\D$,  $g'$ takes the form $g'(z)=\sum_{n=q}^\infty nb_nz^{n-1}$.
Thus, it follows from Corollary \ref{cor3} that
\[
\sum_{ {n=q}}^{\infty} n\left|b_{n}\right| r^{n-1} \leq \sum_{ {n=q}}^{\infty} k n\left|a_{n}\right| r^{n-1} \text { for }
r \leq r_1(q|b_q|/(kq|a_q|))
\]
and integrating this with respect to $r$ gives
\beq \label{equ4}
\sum_{n=q}^{\infty}\left|b_{n}\right| r^{n} \leq k \sum_{n=q}^{\infty}\left|a_{n}\right| r^{n} \text { for }
 r \leq r_1\left(|b_q|/(k|a_q|)\right).
\eeq
In addition, integrating inequality $|g'(z)|^2\leq k^2|h'(z)|^2$ over the circle $|z|=r$, we get (cf. \cite[Lemma 2.1]{KPS})
$$\sum_{n=q}^\infty n^2|b_n|^2r^{2n-2}\leq k^2\sum_{n=q}^\infty n^2|a_n|^2r^{2n-2}\quad
\text{for } r<1,
$$
from which we obtain by integration with respect to $r^2$ that
\beq \label{equ5}
||g||_r^2=\sum_{n=q}^\infty |b_n|^2r^{2n}\leq k^2\sum_{n=q}^\infty |a_n|^2r^{2n}=k^2||h_0||_r^2\quad
\text{for } r<1.
\eeq
Combining \eqref{equ4}, \eqref{equ5} and the inequality in Lemma~D, 
the desired result follows easily.  The remaining part of the proof is easy to obtain by computation (cf. \cite[p.~107]{PVW}).
This completes the proof.
\epf

%
\bthm \label{thm3}
Assume the hypotheses of Lemma  \ref{lem5} and $p\in(0,2]$.
Let $m\in\mathbb{N}$, $|h(0)|=a$, and $q$ be the order of the zero of $h_0$ at $0$.
If $|b_q|\leq1/(2k|a_q|)$, then we have
\beq \label{equ6}
F_{f}^p(z):=|h(z^{m})|^p+E_f(k,r)
\leq 1\quad \text{for }  r=|z|\leq r_{m, k}^{p}(a),
\eeq
where $r_{m, k}^{p}(a)$ is the unique positive root in $(0,1)$ of the equation $\lambda_{m, k}^{p}(a,r)=0$ with
\beq \label{equ7}
\lambda_{m, k}^{p}(a,r)=\{[(1+k)(1-a^2)+1]r-1\}(1+ar^{m})^p+(1-r)(r^m+a)^p.
\eeq
Moreover, for $k=0,~1$ or if $r_{m, k}^p(a)\leq1/3$ for $k\in(0,1)$, 
then the  condition $|b_q|\leq1/(2k|a_q|)$ can be removed and the constant $r_{m, k}^p(a)$ cannot be improved.
\ethm
\bpf
Let us first consider the case of $|b_q|\leq1/(2k|a_q|)$ and $k\neq0$.
Fix $a\in [0,1)$, and observe that the function $\lambda_{m,k}^p(a,r)$ shown in \eqref{equ7} can be rewritten as
\[
\lambda_{m,k}^p(a, r)= (1+a r^{m})^p(1-r)\Lambda_{m,k}^p(a,r),
\]
where
\beq \label{equ8}\Lambda_{m,k}^p(a,r)=\left(\frac{r^{m}+a}{1+a r^{m}}\right)^p+(1+k)(1-a^{2}) \frac{r}{1-r}-1.
\eeq
It is easy to see that  $\Lambda_{m,k}^p(a,r)$ is a strictly increasing function of $r$ in $[0,1)$.
Note that
$$\lambda_{m, k}^{p}(a,0)=a^p-1< 0 ~\mbox{ and }~ \lambda_{m, k}^{p}(a,1)>0.
$$
Clearly, there is a unique positive root $r_{m,k}^p(a)$ in $(0,1)$ of the equation $\lambda_{m, k}^{p}(a,r)=0$.
Further, we have
$$r_{m,k}^p(a)\leq r_{m,k}^2(a) ~\mbox{ for each $p\in(0,2]$},
$$
since $\Lambda_{m,k}^p(a,r)\geq\Lambda_{m,k}^2(a,r)$ for all $r\in[0,1)$.
Simple computation shows that
$$\lambda_{m, k}^{2}(a,r)=(1-a^2)[(1+k)r(1+ar^m)^2-(1-r)(1-r^{2m})],
$$
and thus we have
$$\lambda_{m, k}^{2}(a,1/(2+k))>(1-a^2)[(1+k)r-(1-r)] \big |_{r=1/(2+k)}=0,
$$
which implies $r_{m,k}^2(a)<1/(2+k)$.
If $|b_q|\leq1/(2k|a_q|)$, then we get
$$1/(2+k)\leq1/2\leq r_1(|b_q|/(k|a_q|)).
$$
By Lemmas~E 
and \ref{lem5}, one can obtain that for $|z|=r \leq r_1(|b_q|/(k|a_q|))$,
\beq \label{equ9}
F_{f}^p(z) \leq \left(\frac{r^{m}+a}{1+a r^{m}}\right)^p+(1+k) r \frac{1-a^{2}}{1-r}
=1+\frac{\lambda_{m,k}^p(a, r)}{(1+a r^{m})^p(1-r)},
\eeq
where $\lambda_{m,k}^p(a, r)$ is defined by \eqref{equ7}.
We see that $F_{f}^p(z) \leq 1$ if $\lambda_{m,k}^p(a,r) \leq 0,$ which holds for $r \leq r_{m, k}^p(a),$
where $r_{m,k}^p(a)$ 
is the unique positive root of the equation $\lambda_{m,k}^p(a,r)=0 .$
This proves the inequality \eqref{equ6} if $|b_q|\leq1/(2k|a_q|)$ for $k\neq0$.

For the case of $k=0$, the inequality \eqref{equ6} still holds on the basis of two observations.
One of them is that $g(z)\equiv0$ and so the condition $|b_q|\leq1/(2k|a_q|)$ trivially holds
and thus may be omitted from the theorem.
The other observation is that the inequality \eqref{equ9} holds for $r<1$.

Before checking the sharpness part, we will show a fact, which will be used at later stages as well.
The fact is that  if $r_{m, k}^p(a)\leq1/3$ for $k\neq0$, then the  condition $|b_q|\leq1/(2k|a_q|)$ is not necessary
since $r_1(|b_q|/(k|a_q|))\geq1/3$. Note that $r_{m,1}^p(a)\leq r_{m,1}^2(a)<1/3$ for all $p\in(0,2]$ and all $m\in\mathbb{N}$.
By choosing $f_{a,k}=\omega_a+k\overline{(\omega_a-a)}$ and $z=r$,  we get equality in \eqref{equ9} and thus,
$$F_{f_{a,k}}^p(r)=\left(\frac{r^{m}+a}{1+a r^{m}}\right)^p+(1+k) (1-a^{2}) \frac{r}{1-r}
=1+\frac{\lambda_{m,k}^p(a, r)}{(1+a r^{m})^p(1-r)}.
$$
We see that $F_{f_{a,k}}^p(r) \geq 1$ if and only if $\lambda_{m,k}^p(a,r) \geq 0,$ which holds if and only if $r \geq r_{m, k}^p(a).$
This shows the sharpness  part under the particular conditions in theorem.
This completes the proof of the theorem.
\epf

The condition $r_{m,k}^p(a)\leq1/3$ in Theorem \ref{thm3} is feasible under some simple assumptions,
for instance, $p\in(0,1]$ and $k\geq\frac{1-a}{1+a}$. Indeed, 
 it follows from the proof of Theorem \ref{thm3}
that $r_{m,k}^p(a)\leq r_{m,k}^1(a)$ when $p\in(0,1]$. Direct computations yield
$$\lambda_{m, k}^1(a,r)=(1-a)[(1+k)(1+a)r(1+ar^m)-(1-r)(1-r^m)]
$$
and
$$\lambda_{m, k}^1(a,1/R_k(a))> [(1+k)(1+a)r-(1-r)]\big |_{r=1/R_k(a)}=0,
$$
which implies $r_{m,k}^1(a)<1/R_k(a)$, where $R_k(a) =(1+k)(1+a)+1$.
If $k\geq\frac{1-a}{1+a}$, then we have $r_{m,k}^p(a)\leq r_{m,k}^1(a)\leq1/3$ for $p\in(0,1]$.
Further analysis leads the following result.

\bcor \label{cor4}
Assume the hypotheses of Lemma \ref{lem5} and $p\in(0,2]$.
Let $m\in\mathbb{N}$ and $q$ be the order of the zero of $h_0$ at $0$.
If $|b_q|\leq1/(2k|a_q|)$, then the inequality $F_f^p(z)\leq1$ holds  for $r\leq r_{m,k}^p$,
where  $F_f^p(z)$ is defined by \eqref{equ6}, and $r_{m, k}^p$ is the unique positive root in $(0,1)$ of the equation
 \beq \label{equ10}
\lambda_{m,k}^p(r) =0,
\eeq
where
\be \label{equ10-a}
\lambda_{m,k}^p(r)=2(1+k)r(1+r^{m})-p(1-r)(1-r^{m}).
\ee
Moreover,   for $k=0,~1$ or $p\in(0,1]$, or if $r_{m, k}^p\leq1/3$ for $k\in(0,1)$ and $p\in(1,2]$,
then the condition $|b_q|\leq1/(2k|a_q|)$ can be removed and the constant $r_{m, k}^p$ cannot be improved.
\ecor
\bpf
It is easy to see that $\lambda_{m,k}^p(0)<0$,  $\lambda_{m,k}^p(1)>0$ and $\lambda_{m,k}^p(r)$ is a strictly increasing
function of $r$ in $[0,1)$. Thus, the equation \eqref{equ10} has a unique solution $r_{m,k}^p$ in the interval $(0,1)$.
Simple calculation gives
$$\lambda_{m,k}^p\left(\frac{p}{2(1+k)+p}\right)>[2(1+k)r-p(1-r)]\Big |_{r=\frac{p}{2(1+k)+p}}=0,
$$
which implies
\beq \label{equ11}
r_{m,k}^p\leq\frac{p}{2(1+k)+p}\leq\frac{1}{2+k} \quad  \text{for~all} ~p\in(0,2].
\eeq
If $p\in(0,2]$ and $|b_q|\leq1/(2k|a_q|)$ for $k\neq0$ (resp. $k=0$), it follows from \eqref{equ9} that we have
$F_f^p(z)\leq1+\Lambda_{m,k}^p(a,r)$ for $r\leq1/2\leq r_1(|b_q|/(k|a_q|))$ (resp. $r<1$),
where $\Lambda_{m,k}^p(a,r)$ is defined by \eqref{equ8}.
Now, to show that $F_f^p(z)\leq1$,  it suffices to  prove  the inequality
$\Lambda_{m,k}^p(a,r)\leq0$ for all $a\in[0,1)$, which will be certified  whenever $r\leq r_{m,k}^p$ and $p\in(0,2]$.
This proves the inequality $F_f^p(z)\leq1$ for $r\leq r_{m,k}^p$.

Next for each $r\leq r_{m,k}^p$ and $p\in(0,2]$,
we will prove that $\Lambda_{m,k}^p(a):=\Lambda_{m,k}^p(a,r)$ is an increasing function of $a\in[0,1]$
so   that $\Lambda_{m,k}^p(a)\leq\Lambda_{m,k}^p(1)=0$ for all $a\in[0,1]$.
Elementary calculations provide that
$$(\Lambda_{m,k}^p)'(a)=p(1-r^{2m})\frac{(r^m+a)^{p-1}}{(1+r^ma)^{p+1}}-\frac{2(1+k)ar}{1-r}
$$
and
$$(\Lambda_{m,k}^p)''(a)=p(1-r^{2m})\frac{(r^m+a)^{p-2}}{(1+r^ma)^{p+2}}T_m^p(a,r)-\frac{2(1+k)r}{1-r},
$$
where $$T_m^p(a,r)=(p-1)(1+ar^m)-(p+1)r^m(r^m+a).
$$
Clearly, $(\Lambda_{m,k}^p)''(a)\leq0$ for all $a\in[0,1]$, whenever $p\in(0,1]$.
Hence for $r\leq r_{m,k}^p$,
$$(\Lambda_{m,k}^p)'(a)\geq (\Lambda_{m,k}^p)'(1)=\frac{-\lambda_{m,k}^p(r)}{(1+r^m)(1-r)}\geq0\quad \text{when}~p\in(0,1],
$$
where $\lambda_{m,k}^p(r)$ is defined by \eqref{equ10-a}.
In fact, the assertion $(\Lambda_{m,k}^p)'(a)\geq0$ for $r\leq r_{m,k}^p$ is also   true when $p\in(1,2]$,
which means that $\Lambda_{m,k}^p(a)$ is an increasing function of $a \in[0,1]$ whenever $0<p \leq 2$.
For this, we introduce a function
\[
\Phi(r) 
=\left(\frac{1+r}{1+ar}\right)^2\left(\frac{r+a}{1+ar}\right)^{p-1}:=(\phi_a(r))^2(\varphi_a(r))^{p-1},\quad r \in[0,1).
\]
Simple observations show that no matter $\phi_a(r)$ or $\varphi_a(r)$ for each $a\in[0,1]$, it is an increasing non-negative function of $r$ in $(0,1]$,
so does $\Phi(r)$ when $p>1$. Thus,
$\Phi(r) \geq \Phi(0)=a^{p-1}$ for all $r \in[0,1)$ and for $a \in[0,1]$.
This observation is helpful to derive that for $r \leq r_{m,k}^p$,
\[
\begin{aligned}
(\Lambda_{m,k}^p)'(a) &=p\left(\frac{1-r^m}{1+r^m}\right) \Phi(r^m)-\frac{2 a(1+k) r}{1-r} \\
& \geq a^{p-1}\left[p\left(\frac{1-r^m}{1+r^m}\right)-\frac{2 a^{2-p}(1+k) r}{1-r}\right] \\
& \geq a^{p-1}\left[p\left(\frac{1-r^m}{1+r^m}\right)-\frac{2 (1+k)r}{1-r}\right]=a^{p-1} (\Lambda_{m, k}^p)'(1) \geq 0,
\end{aligned}
\]
since $0 \leq a^{2-p} \leq 1$ for $1<p \leq 2 .$

It remains to show the sharpness part. We choose $f_{a,k}=\omega_a+k \overline{(\omega_a-a)}$ and $z=r$, so we get
$$F_{f_{a,k}}^p(r)=\left(\frac{r^{m}+a}{1+a r^{m}}\right)^p+(1+k) r \frac{1-a^{2}}{1-r}=1+\frac{(1-a)\Psi_{m,k}^p(a, r)}{(1+a r^{m})^p(1-r)},
$$
where
\[
\Psi_{m,k}^p(a, r)=(1-r)(1+ar^m)^p\left[(1+a)\frac{(1+k)r}{1-r}-\frac{1}{1-a}\left(1-\left(\frac{r^m+a}{1+ar^m}\right)^{p}\right)\right].
\]
It is easy to see that $F_{f_{a,k}}^p(r)\geq1$ if and only if $\Psi_{m,k}^p(a,r)\geq 0 .$
In fact, for $r>r_{m,k}^p$ and $a$ close to $1,$ we see that
\[
\lim _{a \rightarrow 1^{-}} \Psi_{m,k}^p(a, r)=(1-r)(1+r^m)^{p}\left[\frac{2 (1+k)r}{1-r}-p\left(\frac{1-r^m}{1+r^m}\right)\right]>0,
\]
which means that the number $r_{m,k}^k$ cannot be improved under particular conditions in the corollary.
Note that $r_{m,k}^p\leq1/3$ for $p\in(0,2]$ if $k\geq p-1$ from \eqref{equ11}.
This finishes the proof of the corollary.
\epf

\brs
Set $|a_0|=a$ and  let $q$ be the order of the zero of $h_0(z)=h(z)-h(0)$ at $0$.

\begin{enumerate}
\item[(1)] The result in Theorem \ref{thm3} (resp. Corollary \ref{cor4}) is still true
if the condition $|b_q|\leq1/(2k|a_q|)$ is replaced by $k=0$ or $r_{m,k}^p(a)\leq r_1(|b_q|/(k|a_q|))$ (resp. $r_{m,k}^p\leq r_1(|b_q|/(k|a_q|)$)
for $k\neq0$. However, the compact expression of $r_{m,k}^p(a)$ (resp. $r_{m,k}^p$) is difficult to state in most cases.

\item[(2)] The result in Theorem \ref{thm3}  still holds for $p>2$ when $k=0$.
This can be seen  from the fact that the function $\Lambda_{m,k}^p$ defined by \eqref{equ8} is also increasing monotonically in $[0,1)$ when $p>2$.
Moreover, the sharpness can be obtained if we choose the function $\omega_a$.
Therefore, Theorem \ref{thm3}  for $k=0$, $m=1$ and $p>0$ coincides with \cite[Lemma~3]{LmsP}, which is a generalization of \cite[Theorem~2]{LLS}.

\item[(3)] Corollary \ref{cor4} for $p=1$ is an improved version of \cite[Theorem~5]{AKP2019}.
Therefore, Theorem \ref{thm3} for $p=1$ is an improved and refined version of \cite[Theorem~5]{AKP2019} under the condition $|b_q|\leq 1/(2k|a_q|)$.
Note that
$$r_{1,0}^p=\frac{p}{\sqrt{4p+1}+p+1}.
$$
Thus, Corollary \ref{cor4} for $k=0$ and $m=1$ leads to \cite[Lemma~2]{LmsP} (i.e. \cite[Lemma~1]{LmsP} with $N=1$),
which is a generalization of \cite[Theorem~1]{LLS} with $N=1$ (an improved version of \cite[Theorem~1]{KP-pre}).
\end{enumerate}
\ers

If we allow $m\rightarrow\infty$ in Theorem \ref{thm3} and Corollary \ref{cor4} in turn,
then we obtain the following two results.

\bcor \label{cor5}
Assume the hypotheses of Lemma  \ref{lem5} and $p\in(0,2]$. Let $q$ be the order of the zero of $h_0$ at $0$,
and $a=|h(0)|$ . If $|b_q|\leq1/(2k|a_q|)$, then the following inequality holds:
\beq \label{equ12}
a^p+E_f(k,r) \leq 1 \quad \text{for}~~ r\leq r_k^p(a):=\frac{1-a^p}{1-a^p+(1+k)(1-a^2)}.
\eeq
Moreover, for $k=0,~1$ or $p\in(0,1]$ {\rm (}resp. for $k\in(0,1)$ and $p\in(1,2)${\rm )},
the condition $|b_q|\leq1/(2k|a_q|)$ can be removed {\rm (}resp. is replaced by $k\geq\frac{1+a^2-2a^p}{1-a^2}${\rm )},
then the above inequality is sharp.
\ecor


\bcor \label{cor6}
Assume the hypotheses of Lemma \ref{lem5} and $p\in(0,2]$.
Let $q$ be the order of the zero of $h_0$ at $0$.
If $|b_q|\leq1/(2k|a_q|)$, then the following inequality holds:
\beq\label{equ13}
|h(0)|^p+E_f(k,r) \leq 1 \quad \text{for}~~  r\leq\frac{p}{2(1+k)+p}.
\eeq
Moreover, for $k=0,~1$ or $p\in(0,1]$ {\rm (}resp. for $k\in(0,1)$ and $p\in(1,2)${\rm )},
the condition $|b_q|\leq1/(2k|a_q|)$ can be removed {\rm (}resp. is replaced by  $k\geq p-1${\rm )}, 
then the above inequality is sharp.
\ecor

It follows from Lemma \ref{lem4} that for $p\in(0,2)$ and $a=|a_0|\in[0,1)$,

$$\inf _{a\in [0,1)} \frac{1+a^2-2a^p}{1-a^2} =p-1 ~\mbox{ and }~
\sup _{a\in [0,1)} \frac{1+a^2-2a^p}{1-a^2}=1.
$$
This means that the condition $k\geq\frac{1+a^2-2a^p}{1-a^2}$ in Corollary \ref{cor5} is reasonable,
and thus 
$r_k^p(a)$ in \eqref{equ12} is no more than $1/3$ under the condition.
Again, it follows from Lemma \ref{lem4} that for $p\in(0,2]$ and $a=|a_0| \in[0,1)$,
$$ \inf _{a\in [0,1)} r_k^p(a)=\inf _{a\in [0,1)}\frac{1}{1+(1+k)t_p(a)}=\frac{p}{2(1+k)+p},
$$
which implies that \eqref{equ13} can be deduced from \eqref{equ12},
where $r_k^p(a)$ is given by \eqref{equ12}.
It is mentioned that the condition $k\geq p-1$ is derived from the inequality $p/[2(1+k)+p]\leq1/3$.

\section{Concluding remarks}
\begin{enumerate}
\item[(1)] Corollaries \ref{cor5} and \ref{cor6} for $k=0$ and $p\in(0,2]$ correspond to
\cite[Remark~1]{PVW2},
which improves Corollary \ref{cor1} and \cite[Proposition~1.4]{bla}.

\item[(2)]  In view of  the second item in the above remarks,
Corollary \ref{cor5} continues to hold for $p>2$ when $k=0$ by applying Lemma \ref{lem5}.
It follows from Lemma \ref{lem4} again that for $p>2$ and $a\in[0,1)$,
\beq\label{equ14}
\inf _{a\in [0,1)} r_k^p(a)=\frac{1}{2+k},
\eeq
where $r_k^p(a)$ is listed in \eqref{equ12}.
Hence, the upper bound of $r$ in Corollary \ref{cor6} is $1/2$ if $k=0$ and $p>2$.

\item[(3)] Note that
$$k=\frac{K-1}{K+1} ~\mbox{ and }~\frac{p}{2(1+k)+p}=\frac{(K+1)p}{(4+p)K+p}.
$$
Thus, Corollary \ref{cor6} for $p=1$ improves \cite[Theorem~1.1]{KPS}.
However, Corollary \ref{cor6} for $p=2$ improves \cite[Theorem~1.2]{KPS} under the condition $|b_q|\leq 1/(2k|a_q|)$.

\item[(4)] Inequality \eqref{equ13} in Corollary \ref{cor6} for $k=1$ and $p=1,~2$ is an improved version of \cite[Corollary~1.4]{KPS}.
Thus, inequality \eqref{equ12} in Corollary \ref{cor5} for $k=1$ and $p=1,~2$ is an improved and refined version of \cite[Corollary~1.4]{KPS}.

\item[(5)]  If $p\in(0,2]$ in \eqref{equ13} is replaced by $p>2$ when $k=1$, then the upper bound of $r$ is $1/3$.
In fact, it follows from Lemma \ref{lem5} that for $r\leq1/3$,
$$|h(0)|^p+E_f(1,r)  \leq  a^p+ (1-a^2)\frac{2r}{1-r}\leq1\quad \text{if}~r\leq r_1^p(a),
$$
where $a=|h(0)|$ and $r_1^p(a)$ is given by \eqref{equ12} with $k=1$.
If we let $k=1$ in \eqref{equ14}, then it is easy to see that the inequality $|h(0)|^p+E_f(1,r)\leq1$ for $r\leq1/3$ when $p>2$.
To see its sharpness, we can consider the function $f(z)=z+\overline{z}$.
A direct computation gives
$$|h(0)|^p+E_f(1,r)=r+r+\frac{1}{1-r}(r^2+r^2)=\frac{2r}{1-r}\geq1
$$
if and only if $r\geq1/3$.
\end{enumerate}

\subsection*{Acknowledgments}
The research of the first author was partly supported by NSFs of China (No. 12071116), the Hunan Provincial Education Department
Outstanding Youth Project (No. 19B079, No. 20B087), the Science and Technology Plan Project of Hunan Province (No. 2016TP1020)
and the Application-Oriented Characterized Disciplines, Double First-Class University Project of Hunan Province (Xiangjiaotong [2018]469).
The  work of the second author is supported by Mathematical Research Impact Centric Support (MATRICS) of
the Department of Science and Technology (DST), India  (MTR/2017/000367).


\end{document}